\documentclass[a4paper, 12pt, openany, fleqn]{amsart}

\usepackage{amssymb}
\usepackage[english]{babel}
\usepackage{blindtext}
\usepackage{graphicx}

\usepackage{hyperref}
\usepackage[a4paper]{geometry}
\usepackage{amsthm}
\usepackage{mathtools}
\usepackage{enumitem}
\usepackage{amsmath}
\usepackage{index}
\usepackage{fancyhdr}
\usepackage{amsbsy}
\usepackage{amsopn}
\usepackage{amsxtra}
\usepackage{amscd}
\usepackage{amstext}
\usepackage{textcomp}
\usepackage[T1]{fontenc}
\usepackage{caption}
\usepackage{float}
\usepackage{inputenc}
\usepackage{tkz-euclide}
\usepackage{mathrsfs}
\usepackage{faktor}
\usepackage{tikz-cd}
\usepackage{tikz}
\usepackage{quiver}
\usetikzlibrary{quotes, babel}
\newtheorem{theorem}{Theorem}[section]
\newtheorem{lemma}[theorem]{Lemma}
\newtheorem{definition}[theorem]{Definition}
\newtheorem{exmp}[theorem]{Example}

\newtheorem*{ack}{Acknowledgments}
\makeindex

\author{Mikhail Malashchuk}
\title{Morphisms of $A_\infty$-algebras and modules}

\date{\today}

\begin{document}

\begin{abstract}
    $A_\infty$-algebra is a chain complex equipped with a coherent multiplicative structure, i.e. a multiplication which is commutative up to (higher) chain homotopy. Similarly $A_\infty$-module over such algebra is equipped with its coherent action. These notions are used to represent cohomological information contained in topological multiplicative structures. One of the key tools of the theory are minimal models, which are useful in cohomological computations. In this thesis author presents certain sufficient conditions for a morphism of $A_\infty$-modules to induce a morphism of their minimal models.
\end{abstract}

\maketitle

\section{Introduction}
$A_\infty$-algebras and modules are widely used tools in homological algebra and geometry. Although they behave nicely over a field and there are many strong results in this case, there are also important constructions over arbitrary rings such as the description of a particular objects in non-commutative geometry by Kontsevich and Soibelman in \cite{MR2596638} or in Heegaard Floer theory in \cite{MR3827056} and \cite{MR4905538}. Therefore the question about generalization of the results over a field seems very natural.\\
The two central results over a field are the following theorems:
    \begin{theorem}[Keller \cite{MR1854636}]
    Let $M$ be an $A_\infty$-module. Then homology $H^*M$ have a structure of an $A_\infty$-module such that:
    \begin{enumerate}
        \item $m_1^{H^*M}\equiv0$, and $m_2^{H^*M}$ is induced by $m_2^M$
        \item There exists a quasi-isomorphism of $A_\infty$-modules $H^*M\longrightarrow M$
    \end{enumerate}
\end{theorem}
\begin{theorem}[Keller \cite{MR1854636}]
    \begin{enumerate}
        \item Homotopy is an equivalence relation on the set of morphisms of $A_\infty$-modules $M\longrightarrow L$
        \item A morphism of $A_\infty$-modules is a quasi-isomorphism if and only if it is a homotopy equivalence
    \end{enumerate}
\end{theorem}
These theorems combined imply functoriality of minimal models of $A_\infty$-module. The goal of this thesis is to present counterexamples to the naive generalizations of these theorems to arbitrary rings and propose some sufficient conditions for which minimal models are functorial. The main result is the following theorem:
\begin{theorem}
    Let $A$ be an $A_\infty$-algebra, $M$ and $L$ be an $A_\infty$-modules over $A$. Assume that $L$ is projective as a grade module and there exists a morphism $f:M\longrightarrow L$ such that $f_1$ is an epimorphism. Then there exists a morphism $\psi :L\longrightarrow M$, which is section of $f$.
\end{theorem}
A dual statement is given in Theorem 3.13 i Section 3 below.\\
We build on the results of \cite{A-infty_rings} to show that without the assumption of projectivity, the result does not hold.\\
The structure of the paper is the following. Section 2 serves as a revision of basic notions of $A_\infty$-algebras, their morphisms and homotopy. We also present known results about them such as stated theorems of Keller and results over an arbitrary ring shown in \cite{A-infty_rings}. Section 3 introduces $A_\infty$-modules in the same way as $A_\infty$-algebras and contains proofs of the main results of the thesis. Section 4 shows how presented theorems imply functoriality of minimal models under certain conditions. The paper is the expanded version of the author's bachelor thesis.\\

\textbf{AI declaration.} The paper was written without using any AI tools.
\begin{ack}
    The author is grateful to Karol Szumiło for the supervision of this thesis and Maciej Borodzik for many fruitful conversations
\end{ack}
\section{$A_\infty$-algebras}
\subsection{Definitions}
We begin by giving the basic definitions of the theory of $A_\infty$-algebras.
\begin{definition}
An $A_\infty$-algebra is a $\mathbb{Z}$-graded $R$-module $A=\bigoplus A^n \text{ for } n\in \mathbb{Z}$ together with homogenous maps $m_n : A^{\otimes n}\longrightarrow A$ of degree $2-n$ for each $n\in\mathbb{N}$, satisfying the following identities:
\begin{align*}
    &\sum_{n=r+s+t} (-1)^{r+st}m_u(\mathbf{1}^{\otimes r}\otimes m_s \otimes \mathbf{1}^{\otimes t}) = 0
\end{align*}
where we set $u=r+1+t$
\end{definition}
\begin{definition}
    A morphism of $A_\infty$-algebras $f:A\longrightarrow B$ is a family $f_n:A^{\otimes n} \longrightarrow B$ of maps of degree $1-n$ such that for each $n\in\mathbb{N}$:
    \begin{align*}
        &\sum_{n=r+s+t} (-1)^{r+st}f_u(\mathbf{1}^{\otimes r}\otimes m_s \otimes\mathbf{1}^{\otimes t})=\\ &=\sum_{1\leq r\leq n\text{ ; } n=i_1 +\dots +i_r} (-1)^s m_r(f_{i_1}\otimes f_{i_2}\otimes\dots\otimes f_{i_r})
    \end{align*}
    where in the first sum we set $u=r+1+t$, and in the second sum we set $s=(r-1)(i_1-1)+(r-2)(i_2-1)+\dots +2(i_{r-2}-1) + (i_{r-1}-1)$
\end{definition}

\begin{definition}
\begin{enumerate}
    \item A morphism $f$ is called a quasi-isomorphism, if $f_1$ is a quasi-isomorphism
    \item A morphism $f$ is call strict, if $f_i \equiv 0$ for all $i\geq 2$
\end{enumerate}
\end{definition}
\begin{definition}
    The composition of morphisms $f:B\longrightarrow C$ and $g:A\longrightarrow B$ is given as
    \begin{align*}
        (f\circ g)_n = \sum_{1\leq r\leq n\text{ ; } n=i_1 +\dots +i_r} (-1)^s f_r \circ (g_{i_1}\otimes g_{i_2}\otimes\dots\otimes g_{i_r})
    \end{align*}
    where we set $s=(r-1)(i_1-1)+(r-2)(i_2-1)+\dots +2(i_{r-2}-1) + (i_{r-1}-1)$
\end{definition}
$A_\infty$-algebras together with such morphisms form a category, in which identity morphisms are strict morphisms such that $f_1$ is an identity. Isomorphisms however don't have to be strict.\\
The following definition introduces an algebraic analogue to topological homotopy and is a natural extension of chain homotopy to the $A_\infty$-algebras.
\begin{definition}
    Let $A$ and $B$ be $A_\infty$-algebras, let $f$ and $g$ be two morphisms $A\longrightarrow B$. A homotopy between $f$ and $g$ is a family of maps $h_i:A^{\otimes i}\longrightarrow B$ of degree $-i$ such that for each $n\in \mathbb{N}$:
    \begin{align*}
        &f_n - g_n =\\ &= \sum (-1)^s m_{r+1+t} (f_{i_1}\otimes\dots\otimes f_{i_r}\otimes h_k\otimes g_{j_1}\otimes\dots\otimes g_{j_t}) +\\ &+ \sum (-1)^{jk+l}h_i(\mathbf{1}^{\otimes j}\otimes m_k \otimes \mathbf{1}^{\otimes l})
    \end{align*}
    where $i_1 +\dots +i_r + k + j_1 +\dots + j_t = n$, $j+k+l=n$ and
    \begin{align*}
        s= t +\sum_{1\leq\alpha\leq t} (1-j_\alpha)(n-\sum_{u\geq \alpha} j_u) + k\sum_{1\leq u \leq r}i_u +\sum_{2\leq\alpha\leq r}(1-i_\alpha)\sum_{u<\alpha}i_u 
    \end{align*}
\end{definition}
\subsection{Known results over a field}
In this section, let $R=\mathbb{K}$ be a field.\\
Over a field are known the following results:
\begin{theorem}[Kadeishvili \cite{MR720689}]
    Let $A$ be an $A_\infty$-algebra. Then homology $H^*A$ have an $A_\infty$-algebra structure such that:
    \begin{enumerate}
        \item $m_1\equiv0$, and $m_2$ is induced by $m_2^A$
        \item There exists a quasi-isomorphism of $A_\infty$-algebras $H^*A\longrightarrow A$
    \end{enumerate}
    Moreover this structure is unique up to an isomorphism of $A_\infty$-algebras
\end{theorem}
\begin{definition}
    The aforementioned structure on $H^*A$ is called a minimal model of an $A_\infty$-algebra.
\end{definition}
\begin{theorem}[Prouté \cite{prout1986}]
    \begin{enumerate}
        \item Homotopy is an equivalence relation on the set of morphisms of $A_\infty$-algebras $A\longrightarrow B$
        \item A morphism of $A_\infty$-algebras is a quasi-isomorphism if and only if it is a homotopy equivalence
    \end{enumerate}
\end{theorem}
These theorems allow to restrict the study only to minimal models of $A_\infty$-algebras and easily pass to a more convenient derived category, which, by Prouté's theorem, is equal to a homotopy category of $A_\infty$-algebras.
\subsection{General results over a ring}
 The results over a field directly use the structure of vector spaces e.g. the direct sum decompositions, which always exist. For modules over arbitrary rings such decomposition does not have to exist and there are known counterexamples for results listed in the previous section.\\
The following example illustrates the failure of Kadeishvili's theorem over a ring without the additional assumptions about the structure of the modules.
\begin{exmp}[Schmid \cite{A-infty_rings}]
    Let $R$ be a domain, in which there exists a non-invertible element $n\in R$. Let $A=R^{2\times 2}$ graded by $A^j = 0$ for $j\in \mathbb{Z}\setminus \{-1,0,1\}$,
    \begin{align*}
        &A^{-1} = \Biggl\{ \begin{pmatrix}
                        0 & 0 \\
                        c & 0
                    \end{pmatrix} \Bigg\vert c\in R\Biggr\}\\
        &A^{0} = \Biggl\{ \begin{pmatrix}
                        a & 0 \\
                        0 & d
                    \end{pmatrix} \Bigg\vert a,d\in R\Biggr\}\\
        &A^{1} = \Biggl\{ \begin{pmatrix}
                        0 & b \\
                        0 & 0
                    \end{pmatrix} \Bigg\vert b\in R\Biggr\}
    \end{align*}
    Then there exists an $A_\infty$-algebra structure on $A$, but there is no compatible minimal $A_\infty$-algebra structure on $H^*A$ quasi-isomorphic in any direction
\end{exmp}
\begin{proof}
    \textbf{Step 1:} There exists an $A_\infty$-algebra structure on $A$:\\
    We set
    \begin{align*}
        &m_1\Biggl(\begin{pmatrix}
            a & b\\
            c & d
        \end{pmatrix} \Biggl) = n\cdot \begin{pmatrix}
            c & d-a\\
            0 & c
        \end{pmatrix}\\
       &m_2\Biggl(\begin{pmatrix}
            a & b\\
            c & d
        \end{pmatrix} \otimes \begin{pmatrix}
            a' & b'\\
            c' & d'
        \end{pmatrix} \Biggl) = \begin{pmatrix}
            a & b\\
            c & d
        \end{pmatrix} \cdot \begin{pmatrix}
            a' & b'\\
            c' & d'
        \end{pmatrix}
    \end{align*}
    The maps $m_1$ and $m_2$ are homogenous of degree $1$ and $0$ accordingly. We also have $m_1^2=0$ and $m_2$ is associative, since it is just a matrix multiplication. Therefore we just have to check the Leibniz rule:
    \begin{flalign*}
        &(m_1\circ m_2)\Biggl(\begin{pmatrix}
            a & b\\
            c & d
        \end{pmatrix} \otimes \begin{pmatrix}
            a' & b'\\
            c' & d'
        \end{pmatrix} \Biggl) = m_1\Biggl(\begin{pmatrix}
            aa'+bc' & ab' + bd'\\
            ca' + dc' & cb' + dd'
        \end{pmatrix}\Biggl)\\
        &= n\cdot \begin{pmatrix}
            ca'+dc' & cb'+dd'-aa'-bc'\\
            0 & ca'+dc'
        \end{pmatrix}\\
        &=n\cdot m_2\Biggl( \begin{pmatrix}
            c & d-a\\
            0 & c
        \end{pmatrix}\otimes\begin{pmatrix}
            a' & b'\\
            c' & d'
        \end{pmatrix} + \begin{pmatrix}
            a & -b\\
            -c & d
        \end{pmatrix}\otimes\begin{pmatrix}
            c' & d'-a'\\
            0 & c'
        \end{pmatrix}\Biggl)\\
        &=m_2\circ (m_1\otimes\mathbf{1} + \mathbf{1}\otimes m_1)\Biggl(\begin{pmatrix}
            a & b\\
            c & d
        \end{pmatrix} \otimes \begin{pmatrix}
            a' & b'\\
            c' & d'
        \end{pmatrix} \Biggl)
    \end{flalign*}
    Thus $(A,m_1,m_2,0,\dots)$ is indeed an $A_\infty$-algebra.\\
    
    \smallskip
    \textbf{Step 2:} Determining the homology $H^*A$:\\
    Since $n$ is not a zero devisor,\\ the set of cycles in $A$ is exactly $Z^*A=\Biggl\{ \begin{pmatrix}
        x & y\\
        0 & x
    \end{pmatrix} \Biggl\vert x,y\in R \Biggl\}$,\\ and the set of boundaries is $B^*A = \Biggl\{n\cdot \begin{pmatrix}
        z & w\\
        0 & z
    \end{pmatrix} \Biggl\vert z,w\in R \Biggl\}$.\\ Hence we have $H^1A\cong H^0A \cong \faktor{R}{nR}$, and the rest is $0$.\\

    \smallskip
    \textbf{Step 3:} There is no minimal $A_\infty$-structure on $H^*A$ together with a quasi-isomorphism $H^*A\longrightarrow A$:\\
    By the choice of $n$ the module $\faktor{R}{nR}$ is a nonzero torsion module annihilated by $n$. $A$, however, is torsion-free, therefore there is no nonzero $R$-linear map $f_1:H^*A\longrightarrow A$, in particular there is no quasi-isomorphism of chain complexes $f_1:(H^*A,0)\longrightarrow (A, m_1)$.\\

    \smallskip
    \textbf{Step 4:} There is no minimal $A_\infty$-structure on $H^*A$ together with a quasi-isomorphism $A\longrightarrow H^*A$:\\
    Suppose that there is such structure $(H^*A,0,m_2',\dots)$ together with a quasi-isomorphism $(f_1,f_2,\dots)$. Then the following equation must be satisfied
    \begin{align*}
        f_1\circ m_2 - f_2\circ(m_1\otimes\mathbf{1}+\mathbf{1}\otimes m_1) = m_2'\circ (f_1\otimes f_1)
    \end{align*}
    We have $m_1 = n\cdot g$, where $g:A\longrightarrow A$, $g\Biggl(\begin{pmatrix}
        a & b\\
        c & d
    \end{pmatrix}\Biggl) = \begin{pmatrix}
        c & d-a\\
        0 & c
    \end{pmatrix}$. Hence we obtain\\ $f_2\circ(m_1\otimes\mathbf{1}+\mathbf{1}\otimes m_1) = n\cdot f_2\circ(g\otimes\mathbf{1}+\mathbf{1}\otimes g) = 0$, since $H^*A$ is annihilated by $n$. Therefore the equation simplifies to
    \begin{align*}
        f_1\circ m_2 = m_2'\circ (f_1\otimes f_1)
    \end{align*}
    Let $\zeta = \begin{pmatrix}
        1 & 0\\
        0 & 1
    \end{pmatrix}$ and $\eta = \begin{pmatrix}
        0 & 1\\
        0 & 0
    \end{pmatrix}$. Then $\zeta$ generates $Z^0A$, a $\eta$ -- $Z^1A$. Since $f_1$ is a quasi-isomorphism, $f_1(\zeta)$ generates $H^0A$, and $f_1(\eta)$ generates $H^1A$, therefore there exists $r\in R$ such that
    \begin{align*}
        -f_1\Biggl(\begin{pmatrix}
        1 & 0\\
        0 & 0
    \end{pmatrix}\Biggl) = r\cdot f_1\Biggl(\begin{pmatrix}
        1 & 0\\
        0 & 1
    \end{pmatrix}\Biggl)
    \end{align*}
    that is $f_1\Biggl(\begin{pmatrix}
        1+r & 0\\
        0 & r
    \end{pmatrix}\Biggl)=0$.\\ Let $\omega = \begin{pmatrix}
        1+r & 0\\
        0 & r
    \end{pmatrix}\otimes\begin{pmatrix}
        0 & 1\\
        0 & 0
    \end{pmatrix} - \begin{pmatrix}
        0 & 1\\
        0 & 0
    \end{pmatrix}\otimes\begin{pmatrix}
        1+r & 0\\
        0 & r
    \end{pmatrix}$.\\ Then we obtain
    \begin{align*}
        &m_2'\circ (f_1\otimes f_1)(\omega)=\\ &= m_2'\Biggl(f_1\Biggl(\begin{pmatrix}
        1+r & 0\\
        0 & r
    \end{pmatrix}\Biggl)\otimes f_1\Biggl(\begin{pmatrix}
        0 & 1\\
        0 & 0
    \end{pmatrix}\Biggl) -\\& -f_1\Biggl(\begin{pmatrix}
        0 & 1\\
        0 & 0
    \end{pmatrix}\Biggl)\otimes f_1\Biggl(\begin{pmatrix}
        1+r & 0\\
        0 & r
    \end{pmatrix}\Biggl)\Biggl) = 0
    \end{align*}
    On the other hand
    \begin{align*}
        (f_1\circ m_2)(\omega) = f_1\Biggl(\begin{pmatrix}
            0 & 1+r\\
            0 & 0
        \end{pmatrix} - \begin{pmatrix}
            0 & r\\
            0 & 0
        \end{pmatrix}\Biggl) = f_1(\eta) \neq 0
    \end{align*}
    since $f_1(\eta)$ generates $H^1A \neq \{0\}$. This is a contradiction.
\end{proof}
An analogue of Kadeishvili's theorem is true if we make additional assumptions about homology of our $A_\infty$-algebra.
\begin{theorem}[Schmid \cite{A-infty_rings}]
    Let $A$ be an $A_\infty$-algebra such that $H^*A$ is a projective $R$-module. Then there is an $A_\infty$-algebra structure on $H^*A$ such that:
    \begin{enumerate}
        \item $m_1\equiv0$, and $m_2$ is induced by $m_2^A$
        \item There exists a quasi-isomorphism of $A_\infty$-algebras $H^*A\longrightarrow A$
    \end{enumerate}
\end{theorem}
\section{$A_\infty$-modules}
\subsection{Definitions}
Let $(A,m_1,m_2,\dots)$ be an $A_\infty$-algebra.
\begin{definition}
    An $A_\infty$-module over $A$ is a $\mathbb{Z}$-graded $R$-module $M$ together with maps
    \begin{align*}
        m_n^M:M\otimes A^{\otimes n-1}\longrightarrow M
    \end{align*}
    of degree $2-n$ for each $n\in\mathbb{N}$ such that
    \begin{align*}
    &\sum_{n=r+s+t} (-1)^{r+st}m_u(\mathbf{1}^{\otimes r}\otimes m_s \otimes \mathbf{1}^{\otimes t}) = 0
\end{align*}
where we set $u=r+1+t$. The term $m_u(\mathbf{1}^{\otimes r}\otimes m_s\otimes\mathbf{1}^{\otimes t})$ is interpreted as $m_u^M(\mathbf{1}^{\otimes r}\otimes m_s\otimes\mathbf{1}^{\otimes t})$ for $r>0$ and as $m_u^M(m_s^M\otimes\mathbf{1}^{\otimes t})$ for $r=0$.
\end{definition}
\begin{definition}
    A morphism of $A_\infty$-modules $f:L\longrightarrow M$ is a family of maps $f_n:L\otimes A^{\otimes n-1}\longrightarrow M$ of degree $1-n$ for each $n\in\mathbb{N}$, satisfying the following identities
    \begin{align*}
        \sum_{n=r+s+t} (-1)^{r+st}f_u\circ (\mathbf{1}^{\otimes r}\otimes m_s\otimes\mathbf{1}^{\otimes t}) = \sum_{n=r+s} (-1)^{(r+1)s}m_u\circ (f_r \otimes \mathbf{1}^{\otimes s})
    \end{align*}
    where in the first sum we set $u=r+1+t$, and in the second sum we set $u=1+s$
\end{definition}
\begin{definition}
    \begin{enumerate}
    \item A morphism $f$ is called a quasi-isomorphism, if $f_1$ is a quasi-isomorphism
    \item A morphism $f$ is called strict, if $f_i \equiv 0$ for all $i\geq 2$
    \item An identity morphism is a strict morphism $f$ such that $f_1$ is an identity
\end{enumerate}
\end{definition}
\begin{definition}
    The composition of morphisms $f:M\longrightarrow N$ and $g:L\longrightarrow M$ is given by
    \begin{align*}
        (f\circ g)_n = \sum_{n=r+s} (-1)^{(r-1)s}f_u\circ(g_r\otimes\mathbf{1}^{\otimes s})
    \end{align*}
    where we set $u=1+s$
\end{definition}
\begin{definition}
    Let $L$ and $M$ be $A_\infty$-modules and let $f$ and $g$ be  morphisms $L\longrightarrow M$. A homotopy between $f$ and $g$ is a family of maps $h_n:L\otimes A^{\otimes n-1}\longrightarrow M$ of degree $-n$ such that for each $n\in\mathbb{N}$
    \begin{align*}
        f_n - g_n = \sum_{n=r+s} (-1)^{rs}m_{1+s}\circ(h_r\otimes\mathbf{1}^{\otimes s}) +\\ + \sum_{n=r+s+t} (-1)^{r+st} h_u\circ(\mathbf{1}^{\otimes r}\otimes m_s\otimes\mathbf{1}^{\otimes t})
    \end{align*}
    where in the second sum we set $u=r+1+t$
\end{definition}
\subsection{Known results over a field}
In this section, let $R=\mathbb{K}$ be a field. We also set an $A_\infty$-algebra $A$.\\
The theorems about $A_\infty$-algebras in this situation are also true for $A_\infty$-modules.
\begin{theorem}[Keller \cite{MR1854636}]
    Let $M$ be an $A_\infty$-module. Then homology $H^*M$ have a structure of an $A_\infty$-module such that:
    \begin{enumerate}
        \item $m_1^{H^*M}\equiv0$, and $m_2^{H^*M}$ is induced by $m_2^M$
        \item There exists a quasi-isomorphism of $A_\infty$-modules $H^*M\longrightarrow M$
    \end{enumerate}
\end{theorem}
The same as for $A_\infty$-algebras such structure is called a minimal model.
\begin{theorem}[Keller \cite{MR1854636}]
    \begin{enumerate}
        \item Homotopy is an equivalence relation on the set of morphisms of $A_\infty$-modules $M\longrightarrow L$
        \item A morphism of $A_\infty$-modules is a quasi-isomorphism if and only if it is a homotopy equivalence
    \end{enumerate}
\end{theorem}
These theorems, analogously to their counterparts for $A_\infty$-algebras, allow us to easily pass to derived category, which has additional structure. An example of the application of this transition is described in detail in \cite{MR1854636}.
\subsection{General results over a ring}
The following examples show that the second theorem also fails over arbitrary ring and emphasize the importance of projectivity/injectivity assumptions for corresponding modules in later theorems.
\begin{exmp}
Let $R=\mathbb{Z}$, $A = \mathbb{Z}$ concentrated in degree $0$, $m_2^A$ is the standard multiplication, and the remaining maps are $0$. Let $M = \dots\rightarrow 0 \rightarrow \mathbb{Z} \xrightarrow{2\cdot} \mathbb{Z}\rightarrow 0 \rightarrow \dots$ together with $m_2^M$ equal to standard multiplication and $m_n^M \equiv 0$ for $n\geq 3$. Then there exists a structure of an $A_\infty$-module on $H^*M$ together with a quasi-isomorphism $f:M\longrightarrow H^*M$, but there is no quasi-isomorphism in the opposite direction
\end{exmp}
\begin{proof}
    \textbf{Step 1:} This structure is an $A_\infty$-module structure on $M$:\\
    The maps $m_1^M$ and $m_2^M$ are of degree $1$ i $0$ accordingly, $(m_1^M)^2=0$, and $m_2^M$ is associative. Thus it suffices to check the Leibniz rule.
    \begin{align*}
        &m_1^M \circ m_2^M(a\otimes b) = m_1^M(ab)=\\ &= 2ab = 2ab + 0= m_2^M(m_1^M(a)\otimes b + a\otimes m_1(b))
    \end{align*}

    \smallskip
    \textbf{Step 2:} Homology $H^*M$:\\
    The homology $H^*M$ is a complex concentrated in degree $0$ with trivial differential and we have $H^*M\cong\mathbb{Z}_2$. It is a $\mathbb{Z}$-module, therefore it is in particular an $A_\infty$-module.\\

    \smallskip
    \textbf{Step 3:} Quasi-isomorphism $M\longrightarrow H^*M$\\
    We set $f = (f_1,0,0,\dots)$, where $f_1$ is the projection $\mathbb{Z}\longrightarrow\mathbb{Z}_2$ in degree $0$ and $0$ everywhere else. It is both a morphism of complexes and a morphism of $\mathbb{Z}$-modules in each degree and since it is exactly the projection on homology, it is a quasi-isomorphism of $A_\infty$-modules.\\

    \smallskip
    \textbf{Step 4:} Lack of quasi-isomorphism $H^*M \longrightarrow M$\\
    $\mathbb{Z}_2$ is a torsion module, and $\mathbb{Z}$ -- torsion-free, therefore there is no nonzero morphism of complexes $H^*M \longrightarrow M$, therefore in particular there is no quasi-isomorphism of $A_\infty$-modules.
\end{proof}

\begin{exmp}
Let $R=\mathbb{Z}$, $A = \mathbb{Z}$ concentrated in degree $0$, $m_2^A$ is a standard multiplication, and the remaining maps are $0$. Let $M = \dots\rightarrow 0 \rightarrow \mathbb{Z} \xrightarrow{\mod{2}} \mathbb{Z}_2\rightarrow 0 \rightarrow \dots$ together with $m_2^M$ equal to a standard multiplication and $m_n^M \equiv 0$ for $n\geq 3$. Then on $H^*M$ there exists an $A_\infty$-module structure together with a quasi-isomorphism $f:H^*M\longrightarrow M$, but there is no quasi-isomorphism in the opposite direction
\end{exmp}
\begin{proof}
    \textbf{Step 1:} This structure is an $A_\infty$-module structure on $M$:\\
     The maps $m_1^M$ and $m_2^M$ are of degree $1$ and $0$ accordingly, $(m_1^M)^2=0$, and $m_2^M$ is associative. Thus it suffices to check the Leibniz rule.
     \begin{align*}
         &m_1^M \circ m_2^M(a\otimes b) = m_1^M(ab) = ab_{\text{mod } 2}=\\ &= (a_{\text{mod } 2})\cdot b + 0 =  m_2^M(m_1^M(a)\otimes b + a\otimes m_1(b))
     \end{align*}

     \smallskip
     \textbf{Step 2:} Homology $H^*M$:\\
     The homology $H^*M$ is a complex concentrated in degree $0$ with trivial differential and we have $H^*M\cong\mathbb{Z}$. It is a $\mathbb{Z}$-module, therefore it is in particular an $A_\infty$-module.\\

     \smallskip
     \textbf{Step 3:} Quasi-isomorphism $H^*M \longrightarrow M$\\
     We $f = (f_1,0,0,\dots)$, where $f_1$ is an inclusion $\mathbb{Z}\longrightarrow 2\mathbb{Z} \hookrightarrow \mathbb{Z}$ in degree $0$ and $0$ everywhere else. It is both a morphism of complexes and morphism of $\mathbb{Z}$-modules in each degree and since it is exactly an inclusion of homology, it is a quasi-isomorphism of $A_\infty$-modules.\\

     \smallskip
     \textbf{Step 4:} Lack of quasi-isomorphism $M\longrightarrow H^*M$\\
     Such quasi-isomorphism has to be a projection on homology. If such projection were to exist $\mathbb{Z}\longrightarrow\mathbb{Z}$, which would be a retraction of an inclusion $\mathbb{Z}\longrightarrow 2\mathbb{Z} \hookrightarrow \mathbb{Z}$, we would obtain $\mathbb{Z}\cong \mathbb{Z}\oplus\mathbb{Z}$, which is a contradiction.
\end{proof}
By making additional assumptions about the structure of modules we can prove the existence of sections or retractions of certain morphisms. The following theorems are dual to each other and each is preceded by a corresponding technical lemma which is analogous to Lemma $1.3.3.3$ in \cite{lefevrehasegawa2003surlesainfinicategories}.
\begin{lemma}
    Let $p:A\longrightarrow D$ such that there exists a section of $p_1$. Then there exists an isomorphism $q:A'\longrightarrow A$ such that $p\circ q$ is a strict morphism
\end{lemma}
\begin{proof}
    Let $A'$ be the same grade module as $A$. We set $q_1 = \mathbf{1}_A$. Assume that we have $q_1,\dots,q_n$ such that
    \begin{align*}
        \sum (-1)^{r-1}p_u(q_r\otimes\mathbf{1}^{\otimes s}) = 0
    \end{align*}
    We now take the $(n+1)$-th identity for compositions without a $p_1(q_{n+1})$ factor and compose it with the section of $p_1$. This is our $q_{n+1}$.\\
    We now want to define the multiplications on $A'$. We set $m_1^{A'}=m_1^A$. Assume that we have $m_1^{A'},\dots,m_n^{A'}$. The $(n+1)$-th identity for $q$ gives us
    \begin{align*}
        q_1(m_{n+1}^{A'}) = -\sum  (-1)^{r+st}q_u(\mathbf{1}^{\otimes r}\otimes m_s \otimes \mathbf{1}^{\otimes t}) +\\ + \sum (-1)^{(r+1)s}m_u^{A}(q_r\otimes\mathbf{1}^{\otimes s})
    \end{align*}
    but $q_1=\mathbf{1}_A$.\\
    We now construct $q^{-1}$. Set $(q^{-1})_1 = \mathbf{1}_A$. Assume that we have $(q^{-1})_1,..,(q^{-1})_n$. From the $(n+1)$-th identity of compositions we obtain
\begin{align*}
    (q^{-1})_{n+1}(q_1\otimes\mathbf{1}^{\otimes n}) = -\sum(-1)^{r-1}(q^{-1})_u(q_r\otimes \mathbf{1}^{s})
\end{align*}
and in the sum there are no factors dependent on $(q^{-1})_{n+1}$, but $q_1 = \mathbf{1}_A$.
\end{proof}
\begin{theorem}
    Let $(A,m_1,m_2,\dots)$ be an $A_\infty$-algebra.\\ Let $(M,m_1^M,m_2^M,\dots)$ and $(L,m_1^L,m_2^L,\dots)$ be $A_\infty$-modules over $A$. Assume that $L$ is projective as a graded module and there exists a morphism $f:M\longrightarrow L$ such that $f_1$ is an epimorphism. Then there exists a morphism $\psi :L\longrightarrow M$, which is section of $f$
\end{theorem}
\begin{proof}
    Let $\psi_1$ be inclusion $L\longrightarrow M$, which exists because of the projectivity of $L$. Then, by the previous lemma, we can without loss of generality assume that $f$ is a strict morphism. From this we get
    \begin{align*}
        f_1(m_n^M) = m_n^L(f_1\otimes\mathbf{1}^{\otimes n-1}) \text{ for each } n\geq 1
    \end{align*}
    and $M = \psi_1(L)\oplus M'$, where $M'=\ker{f_1}$.\\
    Notice that in this situation $\psi = (\psi_1, 0, 0,\dots)$ jest $A_\infty$-morphism. By composing this identity for $f_1$ with $\psi_1$ on the left, and with $\psi_1\otimes\mathbf{1}^{\otimes n-1}$ on the right we obtain
    \begin{align*}
        \psi_1\circ f_1(m_n^M(\psi_1\otimes\mathbf{1}^{\otimes n-1})) = \psi_1\circ m_n^L
    \end{align*}
    but $\psi_1 \circ f_1$ is an identity on $\psi_1(L)$, which completes the proof.
\end{proof}
\begin{lemma}
    Let $j:A\longrightarrow D$ such that there exists a retraction of $j_1$. Then there exists an isomorphism $k:D\longrightarrow D'$ such that $k\circ j$ is a strict morphism
\end{lemma}
\begin{proof} Let $D'$ be the same graded module as $D$. We set $k_1=\mathbf{1}_D$. Assume that we have $k_1,...,k_n$ such that
\begin{align*}
    \sum (-1)^{(r-1)s}k_u(j_r\otimes \mathbf{1}^{\otimes s}) = 0
\end{align*}
We take the $(n+1)$-th identity for compositions without a $k_{n+1}(j_1\otimes \mathbf{1}^{\otimes n})$ factor and compose it with the retraction of $j_1$. This is our $k_{n+1}$.\\
We now want to define the multiplications on $D'$. We set $m_1^{D'} = m_1^{D}$. Assume that we have $m_1^{D'},...,m_n^{D'}$. From the $(n+1)$-th identity for $k$ we obtain
\begin{align*}
    m_{n+1}^{D'}(k_1 \otimes \mathbf{1}^n) = \sum (-1)^{r+st}k_u(\mathbf{1}^{\otimes r}\otimes m_s \otimes \mathbf{1}^{\otimes t}) +\\ + \sum (-1)^{(r+1)s}m_u^{D'}(k_r\otimes\mathbf{1}^{\otimes s})
\end{align*}
but $k_1 = \mathbf{1}_D$.\\
We now construct $k^{-1}$. We set $(k^{-1})_1 = \mathbf{1}_D$. Assume that we have $(k^{-1})_1,..,(k^{-1})_n$. From the $(n+1)$-th identity of compositions we obtain
\begin{align*}
    (k^{-1})_{n+1}(k_1\otimes\mathbf{1}^{\otimes n}) = -\sum(-1)^{r-1}(k^{-1})_u(k_r\otimes \mathbf{1}^{s})
\end{align*}
and in the sum there are no factors dependent on $(k^{-1})_{n+1}$, but $k_1 = \mathbf{1}_D$.
\end{proof}
\begin{theorem}
    Let $(A,m_1,m_2,...)$ be an $A_\infty$-algebra.\\ Let $(M,m_1^m,m_2^m,...)$ and $(L,m_1^L,m_2^L,...)$ be $A_\infty$-modules over $A$. Assume that $L$ is injective as a graded module and there exists a morphism $p:L\longrightarrow M$ such that $p_1$ is a monomorphism. Then there exists a morphism $\varphi :M\longrightarrow L$, which is a retraction of $p$
\end{theorem}
\begin{proof} Let $\varphi_1$ be a projection $M\longrightarrow L$, which exists because of the injectivity of $L$. Therefore we can without loss of generality assume that $p$ is a strict morphism. From this we get:
\begin{align*}
    &p_1(m_n^{L}) = m_n^M(p_1\otimes\mathbf{1}^{\otimes n-1}) \text{ for each } n\geq 1
\end{align*}
and $M = p_1(L)\oplus M'$, where $M' = \ker{\varphi_1}$.\\
Notice that in this situation $\varphi = (\varphi_1,0,0,\dots)$ is $A_\infty$-morphism. Composing this identity with $\varphi_1$ on the left, and with $\varphi_1\otimes\mathbf{1}^{\otimes n-1}$ on the right we obtain, in the same way as in the previous theorem, the desired identity.
\end{proof}

\section{Functoriality of minimal models}
A natural question is, when a morphism $A_\infty$-modules $M\longrightarrow L$ induces a morphism of their minimal models. First of all in this situation one must an make assumption that such models exist. The next natural assumption is the existence of morphisms between this modules and their minimal models. In this situation the theorems from the previous chapter give sufficient conditions for existence of sections or retractions of these morphisms and by taking corresponding compositions we obtain a morphism between the minimal models. This is precisely described in the following theorems:
\begin{theorem}
    Let $M$ and $L$ be $A_\infty$-modules and assume that there is  $A_\infty$-module structures on $H^*M$ and $H^*L$ making them into minimal models together z quasi-isomorphism $H^*M\longrightarrow M$ and $H^*L \longrightarrow L$. Then there exists an induced morphism $H^*M\longrightarrow H^*L$ 
\end{theorem}
\begin{proof}
    From the Theorem $3.13$ we obtain a retraction $L\longrightarrow H^*L$, by composing it with $H^*M\longrightarrow M \longrightarrow L$ we obtain the desired induced morphism.
\end{proof}

The following theorem is dual to the previous and its proof proceeds exactly the same, only instead of the Theorem $3.13$ we apply $3.11$.

\begin{theorem}
    Let $M$ and $L$ be $A_\infty$-modules and assume that there is  $A_\infty$-module structures on $H^*M$ and $H^*L$ making them into minimal models together z quasi-isomorphism $M\longrightarrow H^*M$ and $L \longrightarrow H^*L$. Then there exists an induced morphism $H^*M\longrightarrow H^*L$
\end{theorem}

The following diagram illustrates these two theorems.

\begin{figure}[H]
\[\begin{tikzcd}
	M && L &&& M && L \\
	\\
	{H^*M} && {H^*L} &&& {H^*M} && {H^*L}
	\arrow[from=1-1, to=1-3]
	\arrow["\text{Retraction}", shift left=3, from=1-3, to=3-3]
	\arrow[from=1-6, to=1-8]
	\arrow[from=1-6, to=3-6]
	\arrow[from=1-8, to=3-8]
	\arrow[from=3-1, to=1-1]
	\arrow[dashed, from=3-1, to=3-3]
	\arrow[from=3-3, to=1-3]
	\arrow["{\text{Section}}", shift left=3, from=3-6, to=1-6]
	\arrow[dashed, from=3-6, to=3-8]
\end{tikzcd}\]
    \caption{Application of Theorems 3.11 and 3.13 to inducing morphisms on minimal models}
    \label{fig:1}
\end{figure}
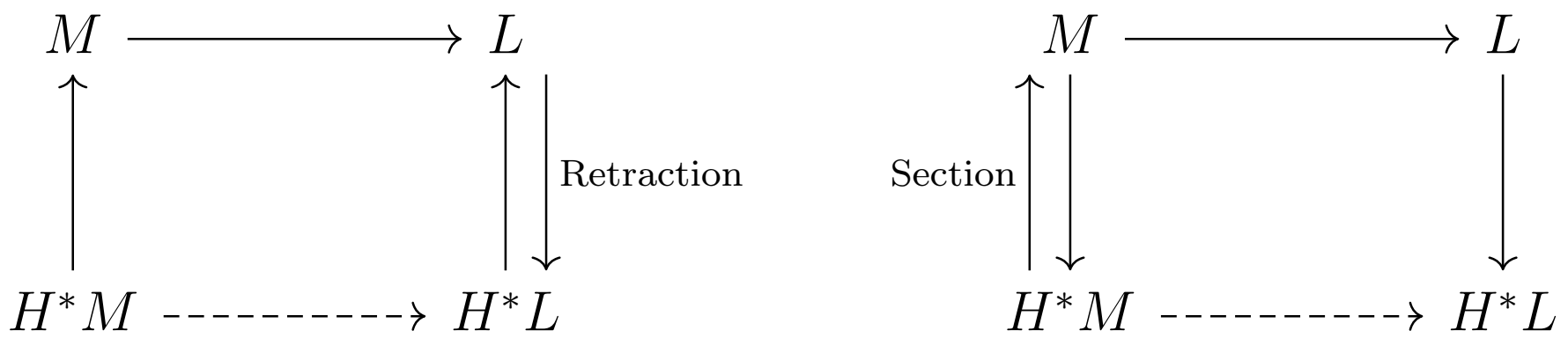

\bibliographystyle{plain}\bibliography{bibliography}
\end{document}